\documentclass[frenchb,english,reqno]{amsart}
\usepackage[T1]{fontenc}
\usepackage[latin1]{inputenc}

\makeatletter

\providecommand{\LyX}{L\kern-.1667em\lower.25em\hbox{Y}\kern-.125emX\@}

 \theoremstyle{definition}
 \newtheorem{defn}{Definition}
 \theoremstyle{plain}    
 \newtheorem{lem}{Lemma} 
 \theoremstyle{remark}
 \newtheorem*{rem*}{Remark}
 \theoremstyle{plain}    
 \newtheorem{thm}{Theorem} 

\newcommand{\textstylehack}{\textstyle}

\newcommand{\NN}{\mathbb {N}}
\newcommand{\TT}{\mathbb {T}}
\newcommand{\RR}{\mathbb {R}}

\newcommand{\EE}{\mathbb {E}}
\newcommand{\VV}{\mathbb {V}}
\newcommand{\PP}{\mathbb {P}}
\newcommand{\one}{\mathbf {1}}
\newcommand{\lebm}{\mathbf {m}}
\newcommand{\brk}{\\ }

\def\essup{\mathop{\operator@font ess\, sup}}
\newcommand{\supp}{\mathrm{supp}\,}

\newcommand{\dist}{\mathrm{dist}\,}
\newcommand{\half}{{\textstyle \frac{1}{2}}}

\newcommand{\quarter}{{\textstyle \frac{1}{4}}}
\newcommand{\threequarters}{{\textstyle \frac{3}{4}}}


\def\endrem{}

\usepackage{babel}
\makeatother
\begin{document}

\title{Random homeomorphisms and Fourier expansions --- the pointwise behavior}

\author{Gady Kozma}

\email{gadykozma@hotmail.com, gadyk@wisdom.weizmann.ac.il}

\curraddr{Gady Kozma\\
The Weizmann Institute of Science\\
Rehovot, Israel.}

\thanks{The research was supported by The Israel Science Foundation (grant
no. 4/01).}

\begin{abstract}
Let $\varphi $ be a Dubins-Freedman random homeomorphism on $[0,1]$
derived from the base measure uniform on $\{x=\frac{1}{2}\}$, and
let $f$ be a periodic function satisfying \[
|f(\delta )-f(0)|=o(\log \log \log \textstylehack \frac{1}{\delta })^{-1}\]
then the Fourier expansion of $f\circ \varphi $ converges at $0$
with probability 1. In the condition on $f$, $o$ cannot be replaced
by $O$.

Also we deduce some 0-1 laws for this kind of problems.
\end{abstract}
\maketitle

\section{Introduction}

\markboth{Random homeomorphisms and Fourier expansions}{The pointwise behavior}
This paper is a continuation of an earlier paper, \cite{KO98}, where
a number of questions related to the Fourier expansions of $f\circ \varphi $
where discussed, most notably conditions under which $S_{n}(f\circ \varphi )$
converges uniformly for a set of $\varphi $'s with probability $1$,
where $S_{n}$ stands for the $n$th Fourier sum. It was proved that
if\[
\omega _{\delta }(f)=o(\log \log \textstylehack \frac{1}{\delta })^{-1}\]
then $S_{n}(f\circ \varphi )$ converges uniformly almost surely,
where $\omega _{\delta }(f)$ stands as usual for the modulus of continuity
of $f$, i.e. \[
\omega _{\delta }(f):=\sup _{|x-y|\leq \delta }|f(x)-f(y)|\quad ,\]
and that this result is sharp (theorems 4 and 6 ibid).

In sections \ref{sect_positive} and \ref{sect_negative} we address
the question of convergence at a specific point. The most obvious
formulation might be {}``under what conditions does $S_{n}(f\circ \varphi )(x)$
converge?'' However, in this formulation it is impossible to get
local conditions on $f$ since $\varphi $ smooths out all the points.
A better formulation uses conditional probability, and reads {}``under
what local conditions on $f$ near $y$ do we have that \[
\left.S_{n}(f\circ \varphi )(x)\, \right|\, \varphi (x)=y\]
converges?'' Essentially, the answer would be the same, i.e.~a triple
log condition, but this formulation incurs a number of technical problems,
so we simplify the proof making use of the fact that $\varphi (0)=0$.
Thus we reached the formulation of the result in the abstract, i.e.

\begin{thm}
\label{thm_logloglog}Suppose $f$ is a continuous function on the
circle satisfying\[
|f(\delta )-f(0)|=o(\log \log \log \textstylehack \frac{1}{\delta })^{-1}\]
Then the Fourier expansion of $f\circ \varphi $ converges at $0$
with probability $1$.
\end{thm}
and this condition is sharp in the following sense: 

\begin{thm}
\label{thm_logloglog_sharp}There exists a continuous function $f$
satisfying\[
|f(\delta )-f(0)|=O(\log \log \log \textstylehack \frac{1}{\delta })^{-1}\]
 for which the Fourier expansion of $f\circ \varphi $ diverges at
$0$ with probability $1$.
\end{thm}
Actually, $f$ may be constructed to satisfy this condition globally,
i.e.\ $\omega _{f}(\delta )=O(\log \log \linebreak [0]\log \frac{1}{\delta })^{-1}$.

It is instructive to contrast these results with the non-stochastic
case. The results of \cite{KO98} are analogues of the Dini-Lipschitz
test \cite[2.71]{17} which gives a sufficient sharp condition for
uniform convergence of $S_{n}(f)$, $\omega _{f}(\delta )=o(\log \frac{1}{\delta })^{-1}$;
for convergence at a point we have the Dini test \cite[2.4]{17} which
gives a sufficient condition $\int \frac{1}{\delta }\omega _{f}(\delta ;x)<\infty $
(again, sharp) where $\omega _{f}(\delta ;x)$ is the modulus of continuity
of $f$ at the point $x$. Thus in the classical case the condition
for pointwise convergence is slightly \textbf{stronger}, or in other
words, a global estimate of $\omega _{f}$ gives better information
about convergence at a specific point than an estimate only at that
point. This behavior, as remarked, does not happen in probabilistic
settings. Of course, we also get a much wider gap, an additional $\log $
factor. We also wish to reiterate remark 4.4i from \cite{KO98}: there
exist functions $f$ satisfying $\omega _{f}(\delta )=O(\log \log \frac{1}{\delta })^{-1}$
such that the Fourier expansion of $f\circ \varphi $ diverges at
a (random) point. This result has no non-probabilistic equivalent.
For a discussion of properties of $S_{n}(f\circ \varphi )$ where
$\varphi $ is non-probabilistic, e.g.~problems such as when $S_{n}(f\circ \varphi )$
might satisfy certain properties for some $\varphi $, all $\varphi $
or a second category set of $\varphi $ see \cite{K83}, \cite{O81}
or \cite{O85}.

Of course, the discussion above does not make much sense without specifying
the probabilistic model for picking $\varphi $, and the group of
homeomorphisms has no Haar measure. We shall be using a model suggested
by Dubins and Freedman \cite{DF65} which uses a base measure $\nu $
on $[0,1]^{2}$. Roughly, a point $(x,y)$ on the graph of $\varphi $
is chosen at random using $\nu $, then this process is repeated for
the rectangles extending from $(0,0)$ to $(x,y)$ and from $(x,y)$
to $(1,1)$ with rescaled versions of $\nu $. repeating this over
and over we get a sequence of points which can, with probability 1,
be closed to a graph of a homeomorphism $[0,1]\rightarrow [0,1]$
with $\varphi (0)=0$ and $\varphi (1)=1$. A proper, though restricted,
definition is provided in section \ref{sub:defrh}. It must be noted,
though, that Dubins and Freedman were not interested in homeomorphisms
but in measures, and considered the Lebesgue-Stieltjes measures $d\varphi $
as random probability measures on $[0,1]$ and studied conditions
under which a typical $d\varphi $ might be singular, atomic and so
on.

Not all Dubins-Freedman measures are born equal, and the most natural
ones are the ones with base measure uniform on $\{x=\frac{1}{2}\}$,
$\{y=\frac{1}{2}\}$ and on $[0,1]^{2}$. See e.g.~\cite{GMW86}
for a specific discussion of these three measures --- they studied
the properties of the set $\varphi (x)=x$ and other interesting facts
about a typical $\varphi $. This paper will be using the first one.
Of course, measures centered on a vertical line are easier to analyse
because one can have an explicit formula for the distributions of
$\varphi (x)$ for dyadic $x$, and sometimes for other $x$'s too,
for example, in the uniform case, $\varphi (\frac{1}{3})$ has the
density function $1-x$ \cite[lemma 1.6]{KO98}. What might be less
clear is that I really need the distribution to be uniform. Indeed,
generalizing the results of \cite{KO98} for measures on $\{x=\frac{1}{2}\}$
which are non-uniform is an open problem. Such a result could be interesting,
for example, in order to play around with the almost-sure H\"older
constant of $\varphi $.

In the last section we discuss the 0-1 law. It turns out that for
this kind of problems, the 0-1 law is not self evident. We shall reduce
the problem to a functional-integral equation (\ref{defp}) which
can be solved by elementary manipulations. This general technique
allows to get 0-1 laws for many problems related to $S_{n}(f\circ \varphi )$:
uniform convergence, pointwise convergence, boundedness of partial
sums etc.

I wish to end this introduction with a question I wasn't even able
to formulate properly. If $I\subset [0,1]$ is a dyadic interval then
the conditional restricted homeomorphism $\psi :=\left.\left.\varphi \right|_{I}\, \right|\, \varphi (\partial I)$
is similar to the original $\varphi $ --- this is the {}``scaling
invariance'', see (\ref{scaling_invariance}) below. If, however,
$I$ is not dyadic then this is no longer true, but $\psi $ still
seems to be very similar to $\varphi $. Many of the results of this
paper and of \cite{KO98} can be reproved for $\psi $. It could be
very interesting (and useful) to prove that for {}``infinitesimal''
problems, $\psi $ and $\varphi $ are equivalent.

\section{Preliminaries}

\subsection{Notations}

We denote by $\TT $ the circle group, which we identify with the
interval $[0,1]$. $\lebm $ denotes the Lebesgue measure on $[0,1]$.
$C$ and $c$ denote absolute positive constants, possibly different,
with $C$ usually pertaining to constants large enough and $c$ to
constants small enough. For a continuous function $f$, $||f||$ denotes
its supremum and $\supp f$ its support.

$\PP $ denotes the probability of some event (with the measure on
the random homeomorphisms defined in the next section). $\EE $ denotes
the expectation of a variable, and $\VV $ its variance. The notation
$X\sim Y$ for two variables means {}``$X$ and $Y$ have the same
distribution''.

Dyadic rational are numbers of the type $k2^{-n}$, $k$ and $n$
integers, and dyadic intervals are intervals of the type $[k2^{-n},(k+1)2^{-n}]$.
For an interval $I:=[a,b]$ the boundary $\partial I$ is the set
$\{a,b\}$. $\left\lfloor x\right\rfloor $ denotes the largest integer
$\leq x$ and $\left\lceil x\right\rceil $ the smallest integer $\geq x$.

$D_{n}$ denotes the Dirichlet kernel on $[0,1]$, i.e.~$\frac{\sin ((2n+1)\pi x)}{\sin (\pi x)}$,
so\[
S_{n}(f;x)=\int _{0}^{1}D_{n}(x-t)\cdot f(t)\, dt\quad .\]

The pointwise modulus of continuity of $f$ at $x$ is defined by
\[
\omega _{f}(x;\delta ):=\sup _{0<|\mu |<\delta }|f(x+\mu )-f(x)|\quad .\]
where for $\delta =0$ we define $\omega _{f}(x;\delta )$ assuming
the function to be periodic.

\subsection{\label{sub:defrh}Random homeomorphisms}

Let's start with the following definition of the particular Dubins-Freedman
measure we will be using, which will be easy to work with. Let $X_{n,k}$
be independent uniform variables in $[0,1]$ for any $n\in \NN $
and any odd $0<k<2^{n}$. We define an increasing function $\varphi $
on the dyadic rational using the following procedure: Start by taking
$\varphi (0)=0$, $\varphi (1)=1$, and\[
\varphi (\half )=X_{1,1}\quad .\]
On the second step, define \[
\varphi (\quarter )=\varphi (\half )\cdot X_{2,1},\quad \varphi (\threequarters )=\varphi (\half )+(1-\varphi (\half ))\cdot X_{2,3}\]
 i.e.~$\varphi (\frac{1}{4})$ and $\varphi (\frac{3}{4})$ are distributed
uniformly on $[0,\varphi (\frac{1}{2})]$ and $[\varphi (\frac{1}{2}),1]$
respectively, and are otherwise independent. We continue this process,
at the $n$th step taking\[
\varphi (k2^{-n}):=\varphi ((k-1)2^{-n})+X_{n,k}\cdot (\varphi ((k+1)2^{-n})-\varphi ((k-1)2^{-n})\quad .\]
This defines $\varphi $ on all dyadic fractions. With probability
$1$, $\varphi $ can be extended to a homeomorphism of $[0,1]$ \cite[theorem 4.1]{DF65}.
We denote this measure by $\PP $, and by $\varphi $ the random change
of variable.

The most useful property of $\varphi $ is {}``scaling invariance'',
which roughly says that for any dyadic interval $I$, $\varphi |_{I}$
behaves like a small copy of $\varphi $. To be more precise,

\begin{lem}
If $I=[k2^{-n},(k+1)2^{-n}]$ is a dyadic interval, then\begin{equation}
\left.\left(\varphi \, |\, \varphi (\partial I)=\{a,b\}\right)\right|_{I}\sim (\varphi \circ L)\cdot (b-a)+a\label{scaling_invariance}\end{equation}
where $L$ is a linear increasing map of $I$ onto $[0,1]$.
\end{lem}
The proof may be found in \cite{GMW86}, theorem 4.6.

Finally, we need the following simple calculation, which can be found
in \cite{KO98} in lemma 1.4 and the remark that follows. For some
constants $K_{1}$ and $K_{2}$ we have\begin{equation}
\PP \left\{ r^{K_{1}}<\varphi (r)<r^{K_{2}}\right\} >1-Cr^{2}\label{phi_is_holder}\end{equation}
for any $r>0$. It will be convenient to assume $K_{2}<1<K_{1}$.

\subsection{And Fourier expansions}

We need the following lemmas, which are deeply related to (though
unfortunately not direct consequences of) theorem 2 from \cite{KO98}:

\begin{lem}
\label{taillemma}For any continuous $f$, $n$, $r>\frac{2}{n}$
and $K>0$,\[
\PP \left(\left|\int _{r}^{1-r}(f\circ \varphi )\cdot D_{n}\right|>K||f||\right)<C\exp \left(-c\frac{\sqrt{nr}}{\log nr}K\right)\quad .\]

\end{lem}
{}

\begin{lem}
\label{headlemma}For any continuous $f$, $n$, interval $I$ and
$K>0$,\[
\PP \left(\left|\int _{I}(f\circ \varphi )\cdot D_{n}\right|>K||f||\right)<Ce^{-e^{cK}}\quad .\]

\end{lem}
and 

\begin{lem}
\label{headlemma_2}For some constant $\beta $, the same $c$ as
above, and any constant $K$,\[
\PP \left(\left|\int _{-r}^{r}(f\circ \varphi )\cdot D_{n}\right|>2K\omega _{f}(0;Cr^{\beta })\right)<Cr^{2}+Ce^{-e^{cK}}\]

\end{lem}
In other words, if theorem 2 of \cite{KO98} gave an estimate of $\int _{0}^{1}(f\circ \varphi )\cdot D_{n}$
then these lemmas give split estimates for the head and the tail.
The proof of lemma \ref{headlemma_2} is an easy corollary to lemma
\ref{headlemma} so let's start with it:

\begin{proof}
Clearly, we may assume $r=2^{-k}$. For each of the segments $[-r,0]$
and $[0,r]$, we use (\ref{phi_is_holder}) ($\beta \equiv K_{2}$),
apply the scaling invariance of $\varphi $ and finally use lemma
\ref{headlemma} for a scaled version of $f$.
\end{proof}
As for the proofs of lemmas \ref{taillemma} and \ref{headlemma},
they follow quite closely the proof of the aforementioned theorem
2, so the rest of this section must be read parallel to it. For lemma
\ref{taillemma}, start from page 1029 ibid. There $||f||=1$ (which
we can also assume here, of course), $I_{k}$ denotes an arc of $\TT $
symmetric around $0$ containing $2k-1$ peaks of the Dirichlet kernel
$D_{n}$ and $Y_{k}:=\int _{I_{k}^{C}}D_{n}\cdot (f\circ \varphi )$.
$I_{k}$ and $Y_{k}$ are connected by the inequality\[
\EE \left(Y_{k}^{2}\: \big |\: \varphi |_{I_{k}}\right)\leq C\frac{\log ^{2}k}{k}\quad ,\]
which is lemma 2.6 ibid. We define $j:=\left\lfloor nr\right\rfloor $
and $\mu :=C\frac{\log j}{\sqrt{j}}$ with $C$ chosen to satisfy\[
\PP \left(|Y_{s}|>\mu \: \big |\: \varphi |_{I_{s}}\right)\leq \textstylehack \frac{1}{4}\quad \forall s\geq j\quad .\]
With this $\mu $ we get

\begin{lem}
\label{lemma_part1_28}If, for a given $\epsilon >0$ and $\nu \geq 1$,
the inequality \[
\PP \left(|Y_{s}|>\nu \mu \: \big |\: \varphi |_{I_{s}}\right)\leq \epsilon \quad \forall s\geq j\]
then\[
\PP \left(|Y_{s}|>(2\nu +2)\mu \: \big |\: \varphi |_{I_{s}}\right)\leq \textstylehack \frac{4}{3}\epsilon ^{2}\]

\end{lem}
The proof is word-for-word identical to the proof of lemma 2.8 ibid.
Now, starting from the definition of $\mu $ we apply lemma \ref{lemma_part1_28}
inductively $l$ times and get that \[
\PP \left(|Y_{s}|>\mu d_{l}\: \big |\: \varphi |_{I_{s}}\right)\leq \left(\textstylehack \frac{4}{3}\right)^{2^{l}-1}\left(\textstylehack \frac{1}{4}\right)^{2^{l}}<\left(\textstylehack \frac{1}{3}\right)^{2^{l}}\]
where the $d_{l}$'s are defined recursively by $d_{1}=1$, $d_{l}=2d_{l-1}+2$.
Clearly $d_{l}<C2^{l}$. Picking a maximal $l$ such that $\mu d_{l}<K$
we get that $2^{l}>cK/\mu $ and lemma \ref{taillemma} follows.\qed 

The proof of lemma \ref{headlemma} is even more similar to that of
theorem 2 from \cite{KO98}, and we shall omit it.

\section{\label{sect_positive}Pointwise convergence}

\emph{Proof of theorem \ref{thm_logloglog}:} Throughout the proof
we shall assume that $f\in C(\TT )$ is some fixed function, that
$||f||\leq 1$ and that $f(0)=0$. We fix $n$ sufficiently large
for the rest of the proof. Define $r=\frac{\log ^{5}n}{n}$. Lemma
\ref{taillemma} will ensure that \begin{equation}
\PP \left(\left|\int _{r}^{1-r}(f\circ \varphi )\cdot D_{n}\right|>\frac{1}{\log n}\right)\leq C\exp \left(-c\frac{\log ^{1.5}n}{\log \log n}\right)<\frac{C}{n^{3}}\label{outside_r}\end{equation}
 Let us now assume that some $m_{1}$ and $m_{2}$ satisfy $m_{1}-m_{2}<\frac{n}{\log ^{6}n}$.
A simple calculation will show\[
|D_{m_{1}}-D_{m_{2}}|=\left|\frac{2\cos ((m_{1}+m_{2}+1)\pi x)\sin ((m_{1}-m_{2})\pi x)}{\sin (\pi x)}\right|<C\frac{n}{\log ^{6}n}\]
so\[
\int _{-r}^{r}|D_{m_{1}}-D_{m_{2}}|<\frac{C}{\log n}\]
 which, combined with (\ref{outside_r}) gives\begin{equation}
\PP \left(\left|\int _{0}^{1}(f\circ \varphi )\cdot (D_{m_{1}}-D_{m_{2}})\right|>\frac{C}{\log n}\right)<\frac{C}{n^{3}}\quad .\label{DM-DN}\end{equation}
Thus, if we only calculate the behavior of $\int (f\circ \varphi )\cdot D_{m}$
on a sequence of $m$'s from $n$ to $2n$ with jumps $\left\lfloor \frac{n}{\log ^{6}n}\right\rfloor $,
we will get a uniform estimate for all $m\in [n,2n]$. Now is the
time to use the log-log-log assumption on $f$. Let $\epsilon (n)\rightarrow 0$
be some sequence converging to $0$ sufficiently slow as to satisfy
\[
\frac{1}{\epsilon (n)}\omega _{f}(0;n^{-\beta })=o(\log \log \log n)^{-1}\]
 Remembering lemma \ref{headlemma_2} (from which we also take the
$\beta $ above), this gives: \[
\PP \left(\left|\int _{-r}^{r}(f\circ \varphi )\cdot D_{n}\right|>\epsilon (n)\right)<Cr^{2}+Ce^{-e^{\Omega (\log \log \log n)}}<\frac{C}{\log ^{8}n}\]
($\Omega $, as usual, denoting the opposite of $o$). We use this
inequality on a sequence of $m$'s which has a length $<C\log ^{6}n$
and throw in (\ref{outside_r}) and (\ref{DM-DN}) to get \[
\PP \left(\exists m\in [n,2n],\; \left|\int _{0}^{1}(f\circ \varphi )\cdot D_{m}\right|>\epsilon (n)+\frac{C}{\log n}\right)<\frac{C}{\log ^{2}n}+n\frac{C}{n^{3}}\]
 and summing these probabilities for $n=2^{k}$ we get the desired
result: that these events happen only for a finite number of $n$'s
for almost every $\varphi $.\qed
\theoremstyle{remark}
\newtheorem*{rems*}{Remarks}
\begin{rems*} 1. Actually, we never used the continuity of $f$. The theorem holds
for any $L^{\infty }$ function satisfying \begin{eqnarray*}
|f(\delta )-f(0^{+})| & = & o(\log \log \log \textstylehack \frac{1}{\delta })^{-1}\\
|f(1-\delta )-f(1^{-})| & = & o(\log \log \log \textstylehack \frac{1}{\delta })^{-1}
\end{eqnarray*}
(for an explanation why $S_{n}(f\circ \varphi )$ is even well defined
for non-continuous $f$, see \cite[lemma 1.3]{KO98}). If $f(0^{+})\neq f(1^{-})$,
we can simply take $f-g$ where $g$ is an appropriate linear function.
$h:=g\circ \varphi $ will be a monotone function, for which we always
have that the Fourier expansion at $x$ converges to $\frac{1}{2}(h^{+}(x)+h^{-}(x))$. 

2. A similar proof shows that for any continuous function $f$, \[
S_{n}(f\circ \varphi ;0)=o(\log \log \log n)\]
and for any $f\in L^{\infty }$,\[
S_{n}(f\circ \varphi ;0)=O(\log \log \log n)\quad .\]
These results too, are sharp.

\end{rems*}

\section{\label{sect_negative}Sharpness.}

This section will be devoted to the proof of theorem \ref{thm_logloglog_sharp}.
Ideologically, the essentials of the proof are contained in the following
heuristics. Examine the following function: \[
f_{n}(t)=\left\{ \begin{array}{ll}
 \sin 2\pi (tn^{k}+\psi _{k}) & t\in [n^{-k},n^{-k+1}],\; 1\leq k\leq e^{n^{4}}\\
 0 & t<n^{-e^{n^{4}}}\end{array}
\right.\]
with some phases $\psi _{k}\in [0,1]$ (usually chosen to make $f_{n}$
continuous). On each interval $\varphi ^{-1}([n^{-k},n^{-k+1}])$,
$f\circ \varphi $ has $n-1$ peaks, and with some small probability
they will be {}``aligned'' with the peaks of some Dirichlet kernel
$D_{r}$. The probability to get a good alignment of $n-1$ peaks
is approximately $e^{-Cn}$, and the variables \[
\varphi |_{\varphi ^{-1}([n^{-k},n^{-k+1}])}\]
 are {}``approximately independent'' so one would expect that in
$>e^{Cn}$ such variables, with big probability this alignment will
happen at least once. In this case, we will have \[
\int _{\varphi ^{-1}([n^{-k},n^{-k+1}])}(f\circ \varphi )\cdot D_{r}\approx c\int |D_{r}|>c\log n\quad .\]
So \[
\sup _{I\subset [0,1],\: r\in \NN }\int _{I}(f\circ \varphi )\cdot D_{r}>c\log n>c\log \log \log n^{e^{n^{4}}}\]
 To make these calculations into a proper proof, we need to do the
following:

\begin{enumerate}
\item \label{explain_align}Explain what it means to {}``get a good alignment
of $f\circ \varphi $ with $D_{r}$'' and calculate the probability.
The calculation will not give $e^{-Cn}$ but a rather weaker estimate
--- hence the element $n^{4}$ in the definition of $f_{n}$.
\item \label{explain_indep}Explain how to overcome the problem that these
{}``approximately independent'' variables are not properly independent.
\item \label{explain_intervals}Explain why it is enough to get a supremum
of $\int _{I}$ for some $I\subset [0,1]$ rather than of $\int _{[0,1]}$.
\item \label{explain_f}Combine the $f_{n}$'s into a single function $f$
which will satisfy the requirements of the theorem. 
\end{enumerate}
We start with issue \ref{explain_intervals}.

\begin{lem}
\label{lemma_I_to_01}Let $||f||\leq 1$, $K>1$, $p$ and $r_{0}<r_{1}$
be given with the condition\[
\PP \left\{ \sup _{\substack{ y\in [0,1]\brk r\in [r_{0},r_{1}]}
}\left|\int _{[0,y]}(f\circ \varphi )\cdot D_{r}\right|>K\right\} >p\]
then\[
\PP \left\{ \sup _{r\in [r_{0},r_{1}]}\left|S_{r}(f\circ \varphi ;0)\right|>\half K\right\} >p-\frac{C}{K}\]

\end{lem}
The proof is practically identical to the proof of lemma 4.5 from
\cite{KO98}, and we shall omit it. The reader might want to skip
to the final steps of the proof of theorem \ref{thm_logloglog_sharp}
on page \pageref{use_lemma_I_to_01} to see how this lemma is used.

To investigate the independence properties of $\varphi $, i.e.~to
explain issue \ref{explain_indep}, let us return to the variables
$X_{n,k}$ defining the measure. For each $i>j\in \NN $ define $\Omega _{i,j}$
to be the $\sigma $-field spanned by \[
\left\{ X_{n,k}\, :\, 2^{-i}<k2^{-n}<2^{-j}\right\} \quad .\]
Clearly, $j>k$ imply that $\Omega _{i,j}$ and $\Omega _{k,l}$ are
independent.

\begin{lem}
\label{lemma_calculus}Let $i>1$ be an integer, $0<\epsilon <1$
and $0<y\leq x\leq 1-\epsilon $. Then\begin{equation}
\PP \left\{ \varphi (1/2)\in [x,x+\epsilon ]\: \big |\: \varphi (2^{-i})=y\right\} >\left(\frac{\epsilon }{2\max (1,|\log y|)}\right)^{i}\label{cond_prob_phi}\end{equation}

\end{lem}
This is a somewhat tedious exercise in calculus. Let us work it out.
A simple calculation (which may be found in \cite{KO98}, (1) page
1022) shows that the distribution function of $\varphi (2^{-i})$
is \begin{equation}
\frac{1}{(i-1)!}\left|\log ^{i-1}y\right|\label{DD_phi}\end{equation}
which, using the scaling invariance of $\varphi $, gives the conditional
distribution function\begin{equation}
\dist \left\{ \varphi (1/2)=x\: \big |\: \varphi (2^{-i})=y\right\} =\left\{ \begin{array}{ll}
 \left|\frac{(i-1)\cdot \log ^{i-2}(x/y)}{x\cdot \log ^{i-1}(y)}\right| & x>y\\
 0 & x\leq y\end{array}
\right.\label{dist_phi_half_dep}\end{equation}
 This distribution (as a function of $x$) is increasing until $x_{0}=ye^{i-2}$
and then decreasing --- which clearly implies that the probability
(\ref{cond_prob_phi}) as a function of $x$ is increasing until some
$x_{1}$ defined by the equality \[
\frac{\log ^{i-2}(x_{1}/y)}{x_{1}}=\frac{\log ^{i-2}((x_{1}+\epsilon )/y)}{x_{1}+\epsilon }\]
and then decreasing, so the minimum is achieved at $x=y$ or, if $x_{1}<1-\epsilon $,
possibly at $x=1-\epsilon $. At $x=y$ we have\[
\PP \left\{ \varphi (1/2)\in [y,y+\epsilon ]\, |\, \varphi (2^{-i})=y\right\} =\frac{\log ^{i-1}\left(1+\frac{\epsilon }{y}\right)}{\log ^{i-1}(y)}\quad .\]
If $\frac{\epsilon }{y}<2$ we estimate $\log \left(1+\frac{\epsilon }{y}\right)>\frac{\epsilon }{2y}>\frac{\epsilon }{2}$
and otherwise $\log \left(1+\frac{\epsilon }{y}\right)>1>\frac{\epsilon }{2}$
so in either case we get (\ref{cond_prob_phi}). At $x=1-\epsilon $,
\[
\PP \left\{ \varphi (1/2)\in [1-\epsilon ,1]\, |\, \varphi (2^{-i})=y\right\} \geq \epsilon \cdot \min _{1-\epsilon \leq x\leq 1}\left|\frac{(i-1)\cdot \log ^{i-2}(x/y)}{x\cdot \log ^{i-1}(y)}\right|\]
and the minimum, remembering $x_{1}<1-\epsilon $, is achieved at
$x=1$ so\[
=\epsilon \cdot \frac{i-1}{|\log y|}\]
and again we get (\ref{cond_prob_phi}).\qed

In the following lemma and its proof the notation $\PP \left(\: \cdot \: \big |\: \varphi (2^{-i})=\varphi _{i}\right)$
is understood as a shorthand for $\lim _{\delta \rightarrow 0}\PP \left(\: \cdot \: \big |\: |\varphi (2^{-i})-\varphi _{i}|<\delta \right)$. 

\begin{lem}
\label{lemma_K}Let $i\in \NN $, let $\varphi _{i}\in [0,e^{-1}]$
and let $\tau $ be an increasing Lipschitz homeomorphism $[2^{-i},1]\rightarrow [\varphi _{i},1]$
with a constant $K$, i.e.~$|\tau (x)-\tau (y)|<K|x-y|$, and let
$0<\epsilon <e^{-1}$. Then\[
\PP \left\{ \left.\max _{2^{-i}\leq x\leq 1}|\varphi (x)-\tau (x)|<\epsilon \right|\varphi (2^{-i})=\varphi _{i}\right\} >c(K)^{1/\epsilon }\cdot \left(\frac{c(K)\epsilon }{|\log \varphi _{i}|}\right)^{C(K)i|\log \epsilon |}\]

\end{lem}
This lemma is a variation on lemma 4.1 from \cite{KO98}, in which
only the $c(K)^{1/\epsilon }$ factor appeared. Think of $c(K)^{1/\epsilon }$
as the {}``main term'', with the other factor meaningful only for
{}``unusual'' cases where $\varphi _{i}$ is very small or $i$
is very large.

\begin{proof}
It is clearly enough to consider $\epsilon =\frac{K+2}{2^{q}}$ where
$q$ is some integer. For any $s\leq q$ we denote\[
A_{s}:=\left\{ \varphi \, :\, |\varphi (j2^{-s})-\tau (j2^{-s})|<2^{-q},\: \forall 2^{s-i}<j<2^{s}\right\} \quad .\]
Let us estimate $\PP (A_{s}\, |\, A_{s-1},\, \varphi (2^{-i})=\varphi _{i})$.
Denote $j^{*}$ to be the minimal $j>2^{s-i}$. If $\varphi \in A_{s-1}$
with some $s\leq q$ then for any odd $j\neq j^{*}$, the probability
of the event\[
|\varphi (j2^{-s})-\tau (j2^{-s})|<2^{-q}\]
 can be estimated from below by\begin{equation}
\frac{2^{-q}}{\left|\varphi \left((j+1)2^{-s}\right)-\varphi \left((j-1)2^{-s}\right)\right|}\geq \frac{2^{-q}}{K2^{-s+1}+2^{-q+1}}\label{j_not_jstar}\end{equation}
and for different $j$'s these events are independent. This estimate
also holds for $j=j^{*}$ if $s>i$, and if $s=i$, $j^{*}$ is even
and therefore irrelevant. Otherwise, for $j=j^{*}$ use lemma \ref{lemma_calculus}
and the scaling invariance of $\varphi $ to get\begin{eqnarray}
\lefteqn{\PP \left\{ \left.|\varphi (j^{*}2^{-s})-\tau (j^{*}2^{-s})|<2^{-q}\, \right|\, \varphi (2^{-i})=\varphi _{i}\right\} } &  & \nonumber \\
 & \qquad  & >\left(\frac{2^{-q}}{2\max \left(1,\left|\log \left(\frac{\varphi _{i}}{\varphi ((j^{*}+1)2^{-s})}\right)\right|\right)}\right)^{i-s}>\left(\frac{2^{-q}}{2|\log \varphi _{i}|}\right)^{i}\quad .\label{j_is_jstar}
\end{eqnarray}
Summing (\ref{j_not_jstar}) and (\ref{j_is_jstar}) (and replacing
$s$ with $s+1$) we get\[
\PP \left(A_{s+1}\, |\, A_{s},\, \varphi (2^{-i})=\varphi _{i}\right)>(K2^{q-s}+2)^{-2^{s}}\cdot \left(\frac{2^{-q-1}}{|\log \varphi _{i}|}\right)^{i}\]
so\begin{eqnarray*}
\lefteqn{\PP \left(A_{q}\, |\, \varphi (2^{-i})=\varphi _{i}\right)>} &  & \\
 & \qquad  & >\left(\frac{2^{-q-1}}{|\log \varphi _{i}|}\right)^{iq}\cdot \prod _{s=0}^{q-1}(K+2)^{-2^{s}}2^{-(q-s)2^{s}}\\
 &  & >\left(\frac{2^{-q-1}}{|\log \varphi _{i}|}\right)^{iq}\cdot \exp \left(-2^{q}\log (K+2)-2^{q}\log 2\sum _{j=1}^{\infty }\frac{j}{2^{j}}\right)\\
 &  & >\left(\frac{c(K)\epsilon }{|\log \varphi _{i}|}\right)^{C(K)i|\log \epsilon |}\cdot \exp \left(-\frac{C(K)}{\epsilon }\right)
\end{eqnarray*}
and clearly $A_{q}$ implies $||\varphi -\tau ||<\epsilon $.
\end{proof}
For the following lemma, we fix $n\in \NN $ large enough and $s_{0}\leq 1$
and inspect the function \[
g(x):=\left\{ \begin{array}{ll}
 f_{n}(xs_{0}^{-1}) & x<s_{0}\\
 0 & \mathrm{otherwise}\end{array}
\right.\]
we define $s_{1}=s_{0}n^{-e^{n^{4}}}$ so that $\supp g=[s_{1},s_{0}]$.

\begin{lem}
\label{lemma_one_I}Let $i>j\in \NN $ and $s_{1}<\varphi _{i}<\varphi _{j}<s_{0}$
satisfy 
\begin{align*} 
4n&<2^{i-j}<n^{K_{1}}\\ 
n^2&<\frac{\varphi_{j}}{\varphi_{i}}<n^{K_{2}} 
\end{align*} 
and let us define the event\[
A_{i,j}:=\left\{ \varphi (2^{-i})=\varphi _{i}\right\} \cap \left\{ \varphi (2^{-j})=\varphi _{j}\right\} \quad .\]
Then\[
\PP \left\{ \left.\exists I\subset [2^{-i},2^{-j}],\, r\in \NN \, :\, \int _{I}(g\circ \varphi )\cdot D_{r}>c\log n\right|A_{i,j}\right\} >e^{-C(K_{1},K_{2})n^{3}}\]

\end{lem}
Note that the above event is in $\Omega _{i,j}$.

\begin{proof}
The conditions on $\varphi _{i}$ and $\varphi _{j}$ imply that for
at least one $k$\begin{equation}
[s_{0}n^{-k},s_{0}n^{-k+1}]\subset [\varphi _{i},\varphi _{j}]\quad .\label{intvl_in_phis}\end{equation}
 Define $k$ to be the least one satisfying (\ref{intvl_in_phis});
$r=(2n)2^{j}$; $\alpha $=$\frac{4}{2r+1}$, $\beta =\frac{2n-4}{2r+1}$,
$I=[\alpha ,\beta ]$; and let us consider the piece-linear homeomorphism
$\tau \, :\, [2^{-i},2^{-j}]\rightarrow [\varphi _{i},\varphi _{j}]$
defined by 
\begin{align*} 
\tau(2^{-i}) &=\varphi_{i}, & \tau(\alpha ) &= (3-\psi_{k})s_{0}n^{-k}, \\ 
\tau(2^{-j}) & = \varphi_{j}, & \tau(\beta ) & = (n-1-\psi_k)s_{0}n^{-k} \quad .
\end{align*} \\
These values were, of course, chosen to ensure $g\circ \tau |_{I}\equiv \sin ((2r+1)\pi x)$.
Simple algebra shows \[
\tau '<C\frac{\varphi _{j}-\varphi _{i}}{2^{-j}-2^{-i}}<C\frac{\varphi _{j}}{2^{-j}}\quad .\]
 Combining this, lemma \ref{lemma_K} and the scaling invariance of
$\varphi $ gives \[
\PP \left\{ \left.\max _{2^{-i}\leq x\leq 2^{-j}}|\varphi (x)-\tau (x)|>\epsilon \right|A_{i,j}\right\} >c^{\varphi _{j}/\epsilon }\cdot \left(\frac{c\epsilon /\varphi _{j}}{|\log (\varphi _{i}/\varphi _{j})|}\right)^{-C(i-j)|\log \epsilon /\varphi _{j}|}\quad .\]
 Taking $\epsilon =n^{-3}\varphi _{j}$ will give \begin{align*} \varphi(\alpha)&>\tau(\alpha)-\epsilon>s_0n^{-k}\left(2-\frac{1}{n}\right) \\ 
\varphi(\beta)&<\tau(\beta)+\epsilon<s_0n^{-k}\left(n-1+\frac{1}{n}\right) 
\end{align*} so\begin{equation}
g\circ \varphi -\sin ((2r+1)\pi x)|_{I}\leq \epsilon \cdot \max _{\varphi I}g'=n^{-3}\varphi _{j}\cdot 2\pi s_{0}^{-1}n^{k}\leq 2\pi n^{-1}\quad .\label{g_phi_is_sin}\end{equation}
The Dirichlet kernel $D_{r}$ and $\sin ((2r+1)\pi x)$ are aligned
in the sense that \[
\int _{I}\sin ((2r+1)\pi x)\cdot D_{r}>c\log n\]
and with (\ref{g_phi_is_sin}), \[
\int _{I}(g\circ \varphi )\cdot D_{r}>c\log n+\frac{C\log n}{n}>c\log n\]
and the probability is\[
>c^{n^{3}}\cdot \left(\frac{cn^{-3}}{K_{2}\log n}\right)^{C\cdot K_{1}\log n\cdot \log n}>e^{-Cn^{3}-C(K_{1},K_{2})\log ^{3}n}>e^{-C(K_{1},K_{2})n^{3}}\qedhere \]

\end{proof}
This lemma is the {}``local'' component of the proof of theorem
\ref{thm_logloglog_sharp}. The complement, the {}``global'' component,
is to show that for typical $\varphi $, many pairs $i$, $j$ satisfying
the conditions above exist. 

\begin{lem}
\label{lemma_cute}Let $0<x\leq y<1$. The probability that $\varphi (2^{-i})\in [0,x]$
where $i$ is the smallest integer satisfying $\varphi (2^{-i})\in [0,y]$
is $\frac{x}{y}$.
\end{lem}
\begin{proof}
Denote this event by $A_{x,y}$. Then\begin{eqnarray*}
\PP A_{x,y} & = & \sum _{i}\PP \left\{ \left(\varphi (2^{-i})\leq x\right)\wedge \left(\varphi (2^{-i+1})>y\right)\right\} \\
 & = & \sum _{i}\int _{y}^{1}\PP \left\{ \left.\varphi (2^{-i})\leq x\, \right|\, \varphi (2^{-i+1})=t\right\} \, d\nu _{i}(t)\\
 & = & \sum _{i}\int _{y}^{1}\frac{x}{t}\, d\nu _{i}(t)\\
 & = & \sum _{i}\int _{y}^{1}\frac{x}{y}\PP \left\{ \left.\varphi (2^{-i})\leq y\, \right|\, \varphi (2^{-i+1})=t\right\} \, d\nu _{i}(t)\\
 & = & \frac{x}{y}\sum _{i}\PP \left\{ \left(\varphi (2^{-i})\leq y\right)\wedge \left(\varphi (2^{-i+1})>y\right)\right\} \\
 & = & \frac{x}{y}\PP \left\{ \exists i\, :\, \left(\varphi (2^{-i})\leq y\right)\wedge \left(\varphi (2^{-i+1})>y\right)\right\} =\frac{x}{y}
\end{eqnarray*}
where the measure $\nu _{i}$ is the distribution of $\varphi (2^{-i+1})$.
\end{proof}
\begin{lem}
\label{lemma_g}For $n$ sufficiently large, for the same $g$ as
above, \[
\PP \left\{ \exists I,\, r\in \NN \, :\, \int _{I}D_{r}\cdot (g\circ \varphi )>c\log n\right\} >1-e^{-n}\]

\end{lem}
\begin{proof}
We use (\ref{phi_is_holder}) for $2^{-d}$ when $d$ is defined by
$d:=\left\lceil \frac{2}{K_{2}}\log _{2}n\right\rceil $ and get \[
\PP \left\{ n^{-K_{3}}<\varphi (2^{-d})<n^{-2}\right\} >1-Cn^{-2}\quad .\]
We need intervals $I_{k}:=[2^{-d(k+1)},2^{-dk}]$ such that $\varphi I_{k}\subset [s_{1},s_{0}]$,
so the first point is to show that many do exist. Lemma \ref{lemma_cute}
ensures that for the random variable $i_{0}$ defined by \[
\left(\varphi (2^{-i_{0}})\leq s_{0}\right)\wedge \left(\varphi (2^{-i_{0}+1})>s_{0}\right)\]
one has\[
\PP \left\{ \varphi (2^{-i_{0}})<e^{-2n}s_{0}\right\} =e^{-2n}\quad .\]
Denote this event by $R_{1}$. Next, define \[
i_{1}:=i_{0}+d\left\lfloor \frac{1}{dK_{1}}\left(\log _{2}n\cdot e^{n^{4}}-2n\log _{2}e\right)\right\rfloor \]
 ($K_{1}$ from (\ref{phi_is_holder})), and using (\ref{phi_is_holder})
and the scaling invariance of $\varphi $ get\[
\PP \left\{ \frac{\varphi (2^{-i_{1}})}{\varphi (2^{-i_{0}})}<n^{-e^{n^{4}}}e^{2n}\right\} <C2^{2(i_{0}-i_{1})}<Ce^{-2n}\quad .\]
Denote this event by $R_{2}$. Between $i_{0}$ and $i_{1}$ we have
$>c_{1}e^{n^{4}}$ intervals $I_{k}$. For each $k$ we define the
event \[
r_{k}:=\neg \left\{ n^{2}<\frac{\varphi (2^{-dk})}{\varphi (2^{-d(k+1)})}<n^{-K_{3}}\right\} \]
so that $\PP r_{k}<Cn^{-2}$ and the $r_{k}$'s are independent. With
these $r_{k}$'s define the variable \[
X:=\#\left\{ k\, :\, I_{k}\subset [i_{0},i_{1}]\wedge r_{k}\right\} \quad .\]
Clearly, $\EE X<c_{1}n^{-2}e^{n^{4}}$ and $\VV X<c_{1}n^{-2}e^{n^{4}}$
so \[
\PP \left\{ X>\half c_{1}e^{n^{4}}\right\} <Ce^{-\half n^{4}}<Ce^{-2n}\quad .\]
Denote this event by $R_{3}$. Finally, we can calculate our probability.
If none of the $R_{i}$'s happen, we have $>\half c_{1}e^{n^{4}}$
intervals $I_{k}$ satisfying the conditions of lemma \ref{lemma_one_I}.
For each $I_{k}$, the behavior of $\left.\varphi |_{I_{k}}\, \right|\, \varphi (\partial I_{k})$
is independent for each $k$ and lemma \ref{lemma_one_I} gives an
estimate of the probability\[
\PP \left\{ \exists I\subset I_{k},r\, :\, \int _{I}D_{r}\cdot (g\circ \varphi )>c\log n\right\} >e^{-Cn^{3}}\]
 ($C$ depends on our $K_{2}$ and $K_{3}$, but is still a constant).
Totally we get\begin{eqnarray*}
\lefteqn{\PP \neg \left\{ \exists I,\, r\in \NN \, :\, \int _{I}D_{r}\cdot (g\circ \varphi )>c\log n\right\} <} &  & \\
 & \qquad \qquad  & (1-e^{-Cn^{3}})^{\left(\half c_{1}e^{n^{4}}\right)}+\PP (R_{1}\cup R_{2}\cup R_{3})<Ce^{-2n}
\end{eqnarray*}
 and the lemma is proved.
\end{proof}

\begin{proof}
[Proof of theorem \ref{thm_logloglog_sharp}] Define values $s_{n}$
and functions $g_{n}$ as follows:\begin{eqnarray*}
g_{n}(x) & := & \left\{ \begin{array}{ll}
 f_{n}(xs_{n}^{-1}) & 4s_{n+1}<x<s_{n}\\
 \mathrm{linear} & x\in [2s_{n+1},4s_{n+1}]\cup [s_{n},2s_{n}]\\
 0 & \mathrm{otherwise}\end{array}
\right.\\
s_{n+1} & := & s_{n}\cdot \quarter n^{-e^{n^{4}}}
\end{eqnarray*}
(take $s_{0}=\frac{1}{4}$) with the relevant $\psi $'s and the linear
portions chosen to make $g_{n}$ continuous. Now pick a sequence $n_{k}\rightarrow \infty $
fast enough as to satisfy, for all $k$,\begin{equation}
\PP \left\{ \exists r\, :\, \left(\left|S_{r}(g_{n_{k}}\circ \varphi ;0)\right|>1\right)\wedge \left(\sum _{l\neq k}\left|S_{r}(g_{n_{l}}\circ \varphi ;0)\right|>\frac{1}{k}\right)\right\} <\frac{1}{k}\label{lll_def_nk}\end{equation}
(this is possible since $S_{r}(g_{n}\circ \varphi ;0)\rightarrow 0$
when $r\rightarrow \infty $ for any fixed $n$ and when $n\rightarrow \infty $
for any fixed $r$). Now define\[
f:=\sum _{k}\frac{1}{\log n_{k}}g_{n_{k}}\]
Clearly, $\omega _{f}(0;\delta )=O(\log \log \log \frac{1}{\delta })^{-1}$.
On the other hand, lemma \ref{lemma_g} ensures that for sufficiently
large $n$,\begin{equation}
\PP \left\{ \exists I\subset [s_{n+1},s_{n}],\, r\in \NN \, :\, \int _{I}D_{r}\cdot (g_{n}\circ \varphi )>c_{1}\log n\right\} >1-e^{-n}\label{prob_one_g}\end{equation}
but $\int _{x}^{y}>c_{1}\log n$ implies that either $\int _{0}^{x}$
or $\int _{0}^{y}>\half c_{1}\log n$. Pick any $r_{1}$ sufficiently
large to allow the restriction $r\in [1,r_{1}]$ in (\ref{prob_one_g}),
and combine this with lemma \label{use_lemma_I_to_01}\ref{lemma_I_to_01}
to get\[
\PP \left\{ \exists r\in \NN \, :\, \left|S_{r}(g_{n}\circ \varphi ;0)\right|>\quarter c_{1}\log n\right\} >1-e^{-n}-\frac{C}{\log n}\quad .\]
and for $n=n_{k}$, again sufficiently large, using (\ref{lll_def_nk})
this event implies \[
\PP \left\{ \sum _{l\neq k}\left|S_{r}(g_{n_{l}};0)\right|>\frac{1}{k}\right\} <\frac{1}{k}\]
 so\[
\PP \left\{ \exists r\in \NN \, :\, \left|S_{r}(f\circ \varphi ;0)\right|>\quarter c_{1}-\textstylehack \frac{1}{k}\right\} >1-e^{-n_{k}}-\frac{C}{\log n_{k}}-\frac{1}{k}\]
and taking $k\rightarrow \infty $ (which clearly forces $r\rightarrow \infty $)
the theorem is done.
\end{proof}
\begin{rem*}
Merely changing the $\frac{1}{\log n_{k}}$ factors in the proof above,
one may get a number of other examples of divergence:\begin{enumerate}
\item[1.] For every $\omega (\delta )=\Omega (\log \log \log \frac{1}{\delta })^{-1}$,
a continuous function $f$ which satisfies $\omega _{f}(\delta ;0)=o(\omega (\delta ))$,
and $S_{n}(f\circ \varphi ;0)$ is almost surely (i.e.~with probability
$1$) unbounded.
\item[2.] For every $\omega (n)=o(\log \log \log n)$, a continuous function
$f$ for which one has $S_{n}(f\circ \varphi ;0)>\omega (n)$ for
infinitely many $n$'s almost surely 
\item[3.] An $L^{\infty }$ function $f$ satisfying $S_{n}(f\circ \varphi ;0)>\log \log \log n$
for infinitely many $n$'s almost surely. \end{enumerate}
\end{rem*}

\section{\label{sect_01}The 0-1 law}

Our aim in this section is to prove claims of the type {}``For any
$f$, the probability that the Fourier expansion of $f\circ \varphi $
converges uniformly (or pointwise, or in 0, or...) is either $0$
or $1$''. As hinted in \cite{KO98} on page 1037, the first step
is to transform the desired property into an {}``interval property'',
for example, to remark that probabilistically, the property 

\[
\sup _{\substack{ n>0\brk t\in [0,1]}
}\left|\int _{[0,1]}D_{n}(t-x)f(x)\, dx\right|<C\]
(i.e.~$f\in U_{0}$, the set of functions with uniformly bounded
Fourier partial sums) is equivalent to \begin{equation}
\sup _{\substack{ I\subset [0,1]\brk n>0\brk t\in [0,1]}
}\left|\int _{I}D_{n}(t-x)f(x)\, dx\right|<C\quad .\label{def_U_tilde}\end{equation}
Denote this set of functions with uniformly bounded {}``interval
Fourier partial sums'' by $\tilde{U}$. That $\PP (f\circ \varphi \in U_{0})=\PP (f\circ \varphi \in \tilde{U})$
was shown in \cite{KO98} in the corollary to lemma 4.5, and $f\in \tilde{U}$
is an interval property, in the following sense:

\begin{defn}
A map $T(f;I)\rightarrow \{0,1\}$ where $f\in C(\TT )$ is a function
and $I\subset [0,1]$ is an interval is called an interval property
if the following conditions hold:
\begin{enumerate}
\item \label{prop_T_Borel}$T$ considered as a map $C(\TT )\times [0,1]^{2}\rightarrow \{0,1\}$
is Borel measurable.
\item $f|_{I}=g|_{I}$ a.e. $\Rightarrow T(f;I)=T(g;I)$.
\item \label{prop_T_multi}$T(f;[x,y])=T(f;[x,t])T(f;[t,y])$ whenever $x\leq t\leq y$;
$T(f;[x,x])=1$
\item \label{prop_T_linear}$T(f\circ L;L^{-1}(I))=T(f;I)$ for any linear
map $L$.
\end{enumerate}
Denote $T(f):=T(f;[0,1])$.
\end{defn}
\endrem When we say that $f\in \tilde{U}$ is an interval property
we mean that the map defined by\[
T(f;I)=1\Leftrightarrow \sup _{\substack{ J\subset I\brk n>0\brk t\in [0,1]}
}\left|\int _{J}D_{n}(t-x)f(x)\, dx\right|<C\]
is an interval property. Property \ref{prop_T_multi} of this $T$
is clear (here the difference between the classes $U_{0}$ and $\tilde{U}$
is crucial). For property \ref{prop_T_linear}, standard arguments%
\footnote{For example, one might show that the difference between the two kernels
(where $n=\left\lfloor \alpha \right\rfloor $) is uniformly bounded.%
} show that the above is equivalent to \[
\sup _{\substack{ J\subset I\brk \alpha \in [0,\infty )\brk t\in \RR }
}\left|\int _{J}\frac{\sin \alpha (t-x)}{t-x}f(x)\, dx\right|<C\]
 for which \ref{prop_T_linear} is clear.

\begin{thm}
\label{theorem_01}If $T(f;I)$ is an interval property and $f$ is
any function, then \[
\PP (T(f\circ \varphi )=1)\in \{0,1\}\]

\end{thm}
\emph{proof:} Let us discuss the following function, defined on $\{0\leq x\leq y\leq 1\}$:\[
p(x,y):=\EE (T(f\circ L_{[x,y]}\circ \varphi ))\]
where $L_{I}$ is the linear increasing mapping of $[0,1]$ onto $I$.
The analysis of $p$ will be based on one equality, (\ref{defp})
below, which we will now prove. \begin{eqnarray*}
p(x,y) & = & \EE \EE \left(T\left(f\circ L_{[x,y]}\circ \varphi ;\left[0,\half \right]\right)T\left(f\circ L_{[x,y]}\circ \varphi ;\left[\half ,1\right]\right)\left|\varphi \left(\half \right)=t\right.\right)\\
 & = & \int _{0}^{1}\EE \left(\left.T\left(f\circ L_{[x,y]}\circ \varphi ;\left[0,\half \right]\right)\right|\varphi \left(\half \right)=t\right)\cdot \\
 &  & \qquad \EE \left(\left.T\left(f\circ L_{[x,y]}\circ \varphi ;\left[\half ,1\right]\right)\right|\varphi \left(\half \right)=t\right)\, dt
\end{eqnarray*}
we now note that, \begin{eqnarray*}
T\left(f\circ L_{[x,y]}\circ \varphi ;\left[0,\half \right]\right) & = & T(f\circ L_{[x,y]}\circ \varphi \circ L_{\left[0,\frac{1}{2}\right]})\\
 & \sim  & T\left(f\circ L_{[x,y]}\circ (t\varphi )\right)\\
 & = & T\left(f\circ L_{[x,x+t(y-x)]}\circ \varphi \right)
\end{eqnarray*}
and similarly\[
T\left(f\circ L_{[x,y]}\circ \varphi ;\left[\half ,1\right]\right)\sim T\left(f\circ L_{[x+t(y-x),y]}\circ \varphi \right)\]
so\[
p(x,y)=\int _{0}^{1}p(x,x+t(y-x))p(x+t(y-x),y)\, dt\]
or, after a change of variable,\begin{equation}
p(x,y)=\frac{1}{y-x}\int _{x}^{y}p(x,t)p(t,y)\, dt\quad .\label{defp}\end{equation}

It might be worth noting that the measurability requirement on $T$
is used only to ensure that $p$ is well defined and measurable on
$[0,1]^{2}$. Thus weaker properties might also do.

First, a technical lemma.

\begin{lem}
\label{simp_h}If $h(x)$ is a bounded function satisfying for every
$x<y_{0}$ \[
h(x)\leq \frac{1}{y_{0}-x}\int _{x}^{y_{0}}h(s)\, ds\]
Then for every $x<t<y_{0}$,\[
h(x)\leq \frac{1}{y_{0}-t}\int _{t}^{y_{0}}h(s)\, ds\]

\end{lem}
\begin{proof}
If not, define \[
s_{0}:=\sup \{s:s<t,\: h(s)\geq h(x)\}\]
and let $s_{n}\rightarrow s_{0}$ be a series satisfying $h(s_{n})\geq h(x)$
(not necessarily different from $s_{0}$). We have for $n$ sufficiently
large,\[
h(x)\leq h(s_{n})\leq \frac{1}{y_{0}-s_{n}}\int _{s_{n}}^{y_{0}}h(s)\, ds\]
so\[
h(x)\leq \frac{1}{y_{0}-t+s_{0}-s_{n}}\left(\int _{s_{n}}^{s_{0}}+\int _{t}^{y_{0}}\right)h(s)\, ds\]
and taking $n\rightarrow \infty $ the lemma is proved.
\end{proof}
\begin{lem}
\label{p_monotone}A measurable function $0\leq p(x,y)\leq 1$ satisfying
\emph{(\ref{defp})} is decreasing in $y$ almost everywhere.
\end{lem}
\begin{proof}
For $x<y<z$, denote\begin{eqnarray*}
\Delta (x,y,z) & := & p(x,z)-p(x,y)\\
\Delta (x,y) & := & \essup _{z\geq y}\Delta (x,y,z)\quad .
\end{eqnarray*}
and define \[
\mu :=\mathrm{ess}\, \sup \Delta (x,y)\]
and assume to the contrary that $\mu >0$. $\Delta $ satisfies the
following:\begin{eqnarray*}
\Delta (x,y,z) & = & \frac{1}{z-x}\left(\int _{x}^{y}p(x,t)\Delta (t,y,z)\, dt+\int _{y}^{z}p(x,t)p(t,z)-p(x,y)\, dt\right)\\
 & \leq  & \frac{1}{z-x}\left(\int _{x}^{y}p(x,t)\Delta (t,y)\, dt+\int _{y}^{z}\Delta (x,y,t)\, dt\right)\quad .
\end{eqnarray*}
We now iterate this inequality. The second iteration looks like that \begin{alignat*}{2} \Delta (x,y,z) & \leq &\;& \frac{1}{z-x}\left( \int _{x}^{y}p(x,t)\Delta (t,y)\, dt\; +\right. \\ \displaybreak[0] &&& \left. \int _{y}^{z}\frac{1}{t-x}\left( \int _{x}^{y}p(x,s)\Delta (s,y)\, ds+\int _{y}^{t}\Delta (x,y,s)\, ds\right) \, dt\right) \leq \\ & = && \frac{1}{z-x}\left( \int _{x}^{y}p(x,t)\Delta (t,y)\left( 1+\int _{y}^{z}\frac{ds}{s-x}\right) \, dt\; +\right. \\ \displaybreak[0] &&& \left. \int _{y}^{z}\Delta (x,y,t)\left( \int _{t}^{z}\frac{ds}{s-x}\right) \, dt\right) = \\ & = && \frac{1}{z-x}\left( \int _{x}^{y}p(x,t)\Delta (t,y)\left( 1+\ln \left( \frac{z-x}{y-x}\right) \right) \, dt\; +\right. \\ &&& \left. \int _{y}^{z}\Delta (x,y,t)\ln \left( \frac{z-x}{t-x}\right) \, dt\right) \end{alignat*}
and, similarly, the $n$th iterate looks like\begin{eqnarray*}
\Delta (x,y,z) & \leq  & \frac{1}{z-x}\left(\int _{x}^{y}p(x,t)\Delta (t,y)\sum _{k=0}^{n-1}\frac{1}{k!}\ln ^{k}\left(\frac{z-x}{y-x}\right)\, dt\; +\right.\\
 &  & \left.\frac{1}{(n-1)!}\int _{y}^{z}\Delta (x,y,t)\ln ^{n-1}\left(\frac{z-x}{t-x}\right)\, dt\right)
\end{eqnarray*}
and, when $n$ tends to infinity, the second term vanishes ($|\Delta |\leq 1$)
and we are left with\[
\Delta (x,y,z)\leq \frac{1}{y-x}\int _{x}^{y}p(x,t)\Delta (t,y)\, dt\]
which is true for all $z$ so\begin{equation}
\Delta (x,y)\leq \frac{1}{y-x}\int _{x}^{y}p(x,t)\Delta (t,y)\, dt\quad .\label{delta}\end{equation}
 Next, fix some small $\epsilon >0$ and get from Lebesgue's density
theorem the existence of a square $[x_{0},x_{0}+\delta ]\times [y_{0},y_{0}+\delta ]$
where $\Delta >\mu -\epsilon $ on a set of measure $>0.9\delta ^{2}$;
and we may also assume that $x_{0}+\delta <y_{0}$ and that $\Delta (x_{0},y_{0}+\delta )>\mu -\epsilon $.
Our contradiction will follow by examining the triangle \[
T:=\{(t,y)\, :\, y_{0}\leq t<y<y_{0}+\delta \}\quad .\]
Now, on one hand we have a set $Y\subset [y_{0},y_{0}+\delta ]$,
$\lebm Y>0.9\delta $ of $y$'s such that for each $y\in Y$ there
exists an $x\in [x_{0},x_{0}+\delta ]$ satisfying $\Delta (x,y)>\mu -\epsilon $
and therefore using lemma \ref{simp_h} for $h(x):=\max \left\{ 0,\, \Delta (x,y)\right\} $
(ignoring, for the moment, the $p$ in inequality (\ref{delta}))
gives\[
\frac{1}{y-y_{0}}\int _{y_{0}}^{y}\max \{0,\Delta (t,y)\}\, dt>\mu -\epsilon \quad \forall y\in Y.\]
This inequality for the average gives a simple measure estimate (assume
$\epsilon <0.1\mu $)\begin{equation}
\lebm \{t\, :\, y_{0}\leq t\leq y,\, \Delta (t,y)>\mu -10\epsilon \}>0.9(y-y_{0})\quad \forall y\in Y\label{defY}\end{equation}
 and on all of $T$ \begin{equation}
\lebm \{(t,y)\in T\, :\, \Delta (t,y)>\mu -10\epsilon \}>0.7\lebm T\quad .\label{delta_is_large}\end{equation}
On the other hand, returning to (\ref{delta}) and inspecting $p$
we get that $\Delta (t,y)>\mu -10\epsilon $ implies\[
\frac{1}{y-t}\int _{t}^{y}p(t,s)\, ds>1-\frac{10\epsilon }{\mu }\]
and as before, \begin{equation}
\lebm \{s\, :\, t\leq s\leq y,\, p(t,s)>1-\frac{100\epsilon }{\mu }\}>0.9(y-t)\label{delta_to_p}\end{equation}
but $y_{0}+\delta \in Y$, which can be combined with (\ref{defY})
and (\ref{delta_to_p}) to get\begin{equation}
\lebm \{(t,s)\in T\, :\, p(t,s)>1-\frac{100\epsilon }{\mu }\}>0.7\lebm T\quad .\label{p_is_large}\end{equation}
Finally, we return to the definition of $\Delta $ and note that $p(x,y)>1-c$
implies $\Delta (x,y)<c$ so we can combine (\ref{delta_is_large})
and (\ref{p_is_large}) to conclude that for some $(x,y)\in T$,\[
\mu -10\epsilon <\Delta (x,y)<\frac{100\epsilon }{\mu }\]
and since $\epsilon $ was arbitrary, the lemma is proved.
\end{proof}
\begin{rem*}
The function \[
p(x,y)=\left\{ \begin{array}{ll}
 0 & y=\frac{1}{2}\\
 1 & \mathrm{otherwise}\end{array}
\right.\]
satisfies (\ref{defp}) but is not monotone everywhere. Thus the {}``almost
everywhere'' in lemma \ref{p_monotone} is not an artifact of the
proof but a property of (\ref{defp}).
\end{rem*}
\begin{lem}
\label{lemma_p_monotone2}A measurable function $0\leq p(x,y)\leq 1$
satisfying \emph{(\ref{defp})} is increasing in $x$ almost everywhere.
\end{lem}
\begin{proof}
Use lemma \ref{p_monotone} for $p'(x,y):=p(1-y,1-x)$.
\end{proof}
We wish to avoid the complexities arising from the fact that $p$
is monotone only almost everywhere. Luckily, all further operations
will be pickings of certain values out of sets of positive measure.
Thus, we can ignore the non-monotone triplets by simply redefining
the notion of picking. Let us call a triplet $x<y<z$ good when $p(x,y)\leq p(x,z)$
and $p(y,z)\leq p(x,z)$; and a triplet $x$, $y$, $z$ is good when
it is good in the right order.

\begin{defn}
\label{def_pick}We say that we pick an $x$ if $x$ satisfies
\begin{enumerate}
\item For almost all $x_{2}$ and $x_{3}$, the triplet $x,x_{2},x_{3}$
is good.
\item If $x_{2}$ has already been picked, then for almost all $x_{3}$,
the triplet $x_{2},x,x_{3}$ is good.
\item If $x_{2}$ and $x_{3}$ have already been picked, then the triplet
$x_{2},x_{3},x$ is good.
\end{enumerate}
\end{defn}
\endrem An induction on lemmas \ref{p_monotone} and \ref{lemma_p_monotone2}
ensures that we can always pick out of every set of positive measure.

\begin{lem}
\label{p_multiplicative}For almost every $x<y<z$, $p(x,z)=p(x,y)p(y,z)$. 
\end{lem}
\begin{proof}
Define \begin{equation}
\Delta (x,y,z):=\left|p(x,z)-p(x,y)p(y,z)\right|\label{def_delta_2}\end{equation}
and assume to the contrary that \begin{equation}
\mu :=\essup \Delta (x,y,z)>0\quad .\label{defmu}\end{equation}
Again, let $\epsilon >0$ be arbitrary, and let \[
[x_{0}-\half \delta ,x_{0}+\half \delta ]\times [y_{0}-\half \delta ,y_{0}+\half \delta ]\times [z_{0}-\half \delta ,z_{0}+\half \delta ]\]
 be a cube where $\Delta >\mu -\epsilon $ on a set of measure $>0.99\delta ^{3}$;
and also assume $x_{0}+\delta <y_{0}$, $y_{0}+\delta <z_{0}$ and
$\Delta (x_{0},y_{0},z_{0})>\mu -\epsilon $. As before, we need a
method to {}``push'' $x$ and $z$ toward $y$. We start from the
simple\begin{equation}
\Delta (x,y,z)\leq \frac{1}{z-x}\left(\int _{x}^{y}p(x,t)\Delta (t,y,z)\, dt+\int _{y}^{z}\Delta (x,y,t)p(t,z)\, dt\right)\quad \label{simpm}\end{equation}
from which we can deduce \theoremstyle{plain} \newtheorem*{sublem}{Sublemma} \begin{sublem}
\label{long_claim}Assume $\Delta (x_{1},y_{1},z_{1})>\mu -\epsilon $
with $|x_{1}-y_{1}|>2\nu $ and $|z_{1}-y_{1}|>2\nu $. Then there
exist $x_{2}$ and $z_{2}$ such that 
\begin{enumerate}
\item $x_{2}<y_{1}<z_{2}$;
\item $\Delta (x_{2},y_{1},z_{2})>\mu -4\epsilon $;
\item $y_{1}-x_{2}<2\nu $;
\item $z_{2}-y_{1}<2\nu $;
\item Either $\nu <y_{1}-x_{2}$ or $\nu <z_{2}-y_{1}$. 
\end{enumerate}
Furthermore, if $x_{1}$, $y_{1}$ and $z_{1}$ are picked in the
sense of definition \ref{def_pick} above then $x_{2}$ and $z_{2}$
are also picked. \end{sublem} \begin{proof}[Subproof] Denoting \[
R_{1}:=\left\{ (x,y_{1},z_{1})\, :\, x\in (x_{1},y_{1})\right\} ,\quad R_{2}:=\left\{ (x_{1},y_{1},z)\, :\, z\in (y_{1},z_{1})\right\} \]
and using (\ref{simpm}) we get\[
\essup _{(x,y,z)\in R_{1}\cup R_{2}}\Delta (x,y,z)>\mu -\epsilon \]
Let us assume that $\essup _{R_{2}}>\mu -\epsilon $. The proof of
the other case will be identical. We denote \[
M(x):=\essup _{z\in [y_{1},z_{1}]}\Delta (x,y_{1},z)\]
(so that $M(x_{1})>\mu -\epsilon $) and using (\ref{simpm}) again,
we have \[
\Delta (x,y_{1},z)\leq \frac{1}{z-x}\left(\int _{x}^{y_{1}}M(t)\, dt+\int _{y_{1}}^{z}M(x)\, dt\right)\]
so\[
M(x)\leq \frac{1}{y_{1}-x}\int _{x}^{y_{1}}M(t)\, dt\]
and we can use lemma \ref{simp_h} for $M$ to get \[
\frac{1}{2\nu }\int _{y_{1}-2\nu }^{y_{1}}M(t)\, dt>\mu -\epsilon \quad .\]
We get a set $X_{1}\subset (y_{1}-2\nu ,y_{1}-\nu )$ of positive
measure with $x\in X_{1}$ satisfying $M(x)>\mu -2\epsilon $, which
implies a set of positive measure $X_{2}\subset (y_{1}-2\nu ,y_{1}-\nu )\times (y_{1},z_{1})$
with $(x,z)\in X_{2}$ satisfying $(x,y_{1},z)>\mu -2\epsilon $.
Let us pick a $z_{3}$ such that $X_{3}:=\left\{ x\, :\, (x,z_{3})\in X_{2}\right\} $
has a positive measure. If $z_{3}<y_{1}+2\nu $, the lemma is proved
--- we denote $z_{2}:=z_{3}$, pick an $x_{2}$ out of $X_{3}$ and
finish. Otherwise, we define\[
M_{2}(z):=\essup _{x\in [y_{1}-2\nu ,y_{1}]}\Delta (x,y_{1},z)\]
and again use lemma \ref{simp_h}, this time for $M_{2}$, to get\[
\int _{y_{1}}^{y_{1}+2\nu }M_{2}(t)\, dt>\mu -2\epsilon \quad .\]
We complete the proof by picking $z_{2}\in (y_{1}+\nu ,y_{1}+2\nu )$
with $M_{2}(z_{2})>\mu -4\epsilon $ and then picking an $x_{2}\in (y_{1}-2\nu ,y_{1})$
satisfying $\Delta (x_{2},y_{1},z_{2})>\mu -4\epsilon $.

\end{proof}

Let us now complete the proof of lemma \ref{p_multiplicative}. First
we use the sublemma for $x_{0}$, $y_{0}$, $z_{0}$ and $\nu =\quarter \delta $.
Denote the resulting values by $x_{1}$ and $z_{1}$. Let us assume
that $y_{0}-x_{1}>\quarter \delta $ --- it will be easy to verify
that the same proof works in the second case. We return to (\ref{simpm})
and observe that $\Delta (x_{1},y_{0},z_{1})>\mu -4\epsilon $ implies
\[
\frac{1}{z_{1}-x_{1}}\left(\int _{x_{1}}^{y_{0}}p(x_{1},t)\, dt+\int _{y_{0}}^{z_{1}}p(t,z_{1})\, dt\right)>1-\frac{4\epsilon }{\mu }\quad .\]
and thus we can pick a $t_{1}\in [\half x_{1}+\half y_{0},y_{0}]$
satisfying $p(x_{1},t_{1})>1-\frac{32\epsilon }{\mu }$. Denote now
$I:=(0.6x_{1}+0.4t_{1},\, 0.4x_{1}+0.6t_{1})$. Since \[
|I|=0.2(t_{1}-x_{1})\geq 0.1(y_{0}-x_{1})>0.025\delta \]
 and since $I\subset [y_{0}-\half \delta ,y_{0}+\half \delta ]$,
we can pick $y_{1}\in I$ such that \[
\lebm \left\{ (x,z)\in [x_{0},x_{0}+\delta ]\times [z_{0},z_{0}+\delta ]\, :\, \Delta (x,y_{1},z)>\mu -\epsilon \right\} >0\quad .\]
 This allows us to proceed and pick $x_{2}$ and $z_{2}$ satisfying
\[
\Delta (x_{2},y_{1},z_{2})>\mu -\epsilon \quad .\]
We use the sublemma again, for $x_{2}$, $y_{1}$, $z_{2}$ and $\nu =0.05\delta $.
Denoting the output of the claim by $x_{3}$ and $z_{3}$ we are finally
faced with the following situation: \begin{gather*} 
x_1<x_3<y_1<z_3<t_1<y_0<z_1 \\ 
p(x_1,t_1)>1-\frac{32\epsilon}{\mu}\\ 
\Delta(x_3,y_1,z_3)>\mu-4\epsilon\quad.
\end{gather*} 
This, however, is a contradiction to the assumption $\mu >0$ since
the monotonicity of $p$ gives\[
p(x_{3},y_{1}),\, p(y_{1},z_{3})>1-\frac{32\epsilon }{\mu }\]
so\[
\Delta (x_{3},y_{1},z_{3})<1-\left(1-\frac{32\epsilon }{\mu }\right)^{2}<\frac{64\epsilon }{\mu }\]
and, since $\epsilon $ was arbitrary, $\mu $ must be zero. 
\end{proof}
The fact that $p(x,y)$ is multiplicative only almost everywhere requires
us to use a variation on the standard 0-1 law. The formulation follows:

\begin{lem}
\label{lemma_regular_01}Let $\Omega =\prod \Omega _{n}$ be a (product)
probability space and $X$ a random variable defined on $\Omega $
such that for \emph{almost} every $\omega _{1},\omega _{1}'\in \Omega _{1},\ldots ,\omega _{n},\omega _{n}'\in \Omega _{n}$,\[
\EE (X|\omega _{1},\ldots ,\omega _{n})=\EE (X|\omega _{1}',\ldots ,\omega _{n}')\]
then $X$ is similar to a constant. In particular, if $X=\one _{A}$,
then $\PP (A)=0$ or $\PP (A)=1$.
\end{lem}
The proof is identical to the proof of the standard 0-1 law --- see
e.g.~\cite[page 7]{K}.

\emph{Proof of theorem \ref{theorem_01}:} We want to use lemma \ref{lemma_regular_01}
with the independent variables $X_{n,k}$. Clearly, we may assume
that the number of variables in the lemma is $2^{N}-1$. Now, taking
$\EE (\, \cdot \, |\left\{ X_{n,k}=\omega _{n,k}\right\} )$ for $n=1,\ldots ,N$
and $1<k<2^{n}$ is identical to taking \[
\EE \left(\, \cdot \, \left|\left\{ \varphi (k2^{-N})=s_{k}\right\} _{k=0}^{2^{N}}\right.\right)\]
(write $s_{0}=0$ and $s_{2^{N}}=1$) and then\begin{eqnarray*}
\lefteqn{\EE \left(T(f\circ \varphi ;[0,1])\left|\left\{ \varphi (k2^{-N})=s_{k}\right\} _{k=0}^{2^{N}}\right.\right)=} &  & \\
 & \qquad \qquad  & =\EE \left(\left.\prod _{k=0}^{2^{N}-1}T(f\circ \varphi ;[k2^{-N},(k+1)2^{-N}])\right|\left\{ \varphi (k2^{-N})=s_{k}\right\} _{k=0}^{2^{N}}\right)\\
 &  & =\prod _{k}\EE \left(T(f\circ \varphi ;[k2^{-N},(k+1)2^{-N}])\left|\left\{ \varphi (l2^{-N})=s_{l}\right\} _{l=k,k+1}\right.\right)\\
 &  & =\prod _{k}\EE \left(T(f\circ L_{[s_{k},s_{k+1}]}\circ \varphi )\right)\\
 &  & =\prod _{k}p(s_{k},s_{k+1})=p(0,1)\quad \mathrm{for}\, \mathrm{a}.\mathrm{e}.\, \left\{ s_{k}\right\} 
\end{eqnarray*}
 and the theorem is proved.\qed 

Theorem \ref{theorem_01} can be applied to a number of harmonic properties
of $f\circ \varphi $. Let us name a few, without proofs:

\begin{itemize}
\item Uniform convergence of $S_{n}(f\circ \varphi )\rightarrow f\circ \varphi $.
One possible corresponding interval property is \[
T(f;I)=1\Leftrightarrow \forall J\subset I,\lim _{n\rightarrow \infty }\int _{J}D_{n}(x-t)\cdot \left(f(\varphi (t))-f(\varphi (x))\right)\, dt=0\]
uniformly in $x$. Showing that this is probabilistically equivalent
to $f\in U(\TT )$ is similar to the proof that $U_{0}$ is equivalent
to $\tilde{U}$.
\item Pointwise divergence on an infinite/uncountable/dense/second\linebreak[0]
category set. All these properties (or their complements) are interval
properties to begin with, so theorem \ref{theorem_01} applies directly.
\item Pointwise convergence everywhere.
\item Pointwise bounded Fourier partial sums.
\item For any $\psi (n)\nearrow \infty $, \[
S_{n}(f\circ \psi )=o(\psi (n))\]
uniformly or pointwise everywhere.
\end{itemize}

\end{document}